\documentclass[a4paper,11pt]{amsart}
\usepackage{hyperref,latexsym}
\usepackage{enumerate}

\theoremstyle{plain}
\newtheorem{theorem}{Theorem}[section]
\newtheorem{lemma}[theorem]{Lemma}

\theoremstyle{definition}
\newtheorem{definition}[theorem]{Definition}

\theoremstyle{remark}
\newtheorem{remark}{Remark}

\begin{document}

\title[Richness]
      {Richness or Semi-Hamiltonicity  of  quasi-linear systems which are not in evolution form}

\date{23 January 2011}
\author{Misha Bialy}
\address{School of Mathematical Sciences, Raymond and Beverly Sackler Faculty of Exact Sciences, Tel Aviv University,
Israel} \email{bialy@post.tau.ac.il}
\thanks{Partially supported by ISF grant 128/10}

\subjclass[2000]{35L65,35L67,70H06 } \keywords{Rich,
Semi-Hamiltonian, genuine nonlinearity, blow-up, Systems of
Hydrodynamic type}

\begin{abstract}

 The aim of this paper is to consider quasi-linear systems which are not
in the form of evolution equations. We propose new condition of
Richness or Semi-Hamiltonicity for such a system and prove that the
blow up analysis along characteristic curves can be performed for it
in an analogous manner. This opens a possibility to use this ansatz
also for geometric problems. We apply the results to the problem of
Polynomial integral for geodesic flows on the 2-torus.
\end{abstract}

\maketitle

%%%%%%%%%%%%%%%%%%%%%%%%%%%%%%%%%%%%%%%%%%%%%%%%%%%%%%%%%%%%%%%%%%%%%%%%%%
\section{Motivation and the result}
\label{sec:intro} Consider a quasi-linear system for vector function
$u(x,y)=(u_{1},...,u_{n})$ which has the following form
$$A(u)u_{x}+B(u)u_{y}=0 \eqno {(1)}$$
It may happen in practice that one of the matrices $A(u)$ and $B(u)$
can degenerate somewhere (and even both of them can degenerate
somewhere).

 Our main assumption about these matrices is that the homogeneous polynomial $P$
in $\alpha,\beta$ is not a zero polynomial at any point $(x,y):$
$$P=\det(\alpha B -\beta A),\ deg(P)=n. \eqno (P)
$$
This assumption is obviously satisfied if one of the matrices $A(u)$
and $B(u)$ is non-degenerate, however we shall assume everywhere the
weaker version-$(P)$. We shall see in the example that $(P)$ is in
fact the correct assumption.

Moreover, we shall assume in the following that the system is
strictly hyperbolic that is the polynomial $P$ has $n$ distinct
roots $[\beta_i : \alpha_i].$ We define unite characteristic vector
fields on the plane $\mathbf{R}^2(x,y)$ by
$$ v_i=\cos \phi_i \partial_x+\sin \phi_i \partial_y,$$
where the angles $\phi_i$, we shall call them characteristic angles,
are such that $[\sin \phi_i:\cos \phi_i]= [\beta_i:\alpha_i].$

Notice that if one of the matrices, say $A$ is non-degenerate then
the system is equivalent to one in the evolution form. The notion of
the system to be  Rich or Semi-Hamiltonian (see
\cite{Serre},\cite{Sev} and \cite{Tsarev},\cite{DN}, we shall call
them Rich for the sake of brevity) for evolution system says that it
can be written in Riemann invariants (diagonal form)
$$(r_i)_x +\lambda_i(r_1,...,r_n)(r_i)_y=0, i=1,...,n,$$  and moreover the
eigenvalues $\lambda_i= \beta_i/\alpha_i$ of $A^{-1}B$ satisfy the
following identities:
$$\partial_{r_k}\left(\frac{\partial_{r_i}\lambda_j}{\lambda_i-\lambda_j}\right)=\partial_{r_i}
\left(\frac{\partial_{r_k}\lambda_j}{\lambda_k-\lambda_j}\right).
\eqno(R)$$ This condition allows one to perform blow-up analysis
along characteristics as it is shown in \cite{Serre} and used for
example in \cite{B1}. It was proved by B. Sevennec \cite{Sev} and
later understood in differential-geometric terms \cite{F} that
strictly Hyperbolic system in evolution form which is written in
Riemann invariants is Rich, if and only if there are local
coordinates in which the system takes the form of conservation laws.

The unsatisfactory thing, however,with the condition (R) is the fact
that for the system (1) characteristic curves can pass from the
chart where $A$ is non-degenerate to the chart where $B$ is
non-degenerate or even reach those points where both matrices
degenerate. This does-not allow to use the analysis of corresponding
Riccati equations for all times.

We propose the following generalization of the Richness condition
whose naturality  we shall justify below:

\begin{definition}
We call system (1) Rich if it can be written in the diagonal
form
$$L_{v_i}r_i=\cos \phi_i (r_i)_x+ \sin \phi_i (r_i)_y=0,\ i=1,...,n \eqno(2)$$ for a
regular change of variables $(u_{1},...,u_{n}) \rightarrow
(r_{1},...,r_{n})$ and the following conditions on the
characteristic angles $\phi_i(r_{1},...,r_{n})$ holds true
$$\partial_{r_k}\left(\frac{\partial_{r_i}\phi_j}{\tan(\phi_i-\phi_j)}\right)=\partial_{r_i}
\left(\frac{\partial_{r_k}\phi_j}{\tan(\phi_k-\phi_j)}\right)
\eqno(\Phi)$$
 \end{definition}
\smallskip
 It is important fact that this definition is invariant with respect to rotations
 of the plane. We shall continue to call $r_i$ in (2) by Riemann
 invariants. Our first result is the following

\begin{theorem} If system (1) satisfying (P) is Rich according Definition
1.1, so that the conditions (2),($\Phi$) hold true, then the
derivatives of i-th Riemann invariant $w_i=L_{v_i^\bot}r_i$ in the
orthogonal direction to characteristics satisfy the following
Riccati equation:
$$L_{v_i}(\exp{(-G_i)}w_i)+ \exp{(G_i)}\partial_{r_i}(\phi
_i)(\exp({-G_i})w_i)^2=0,$$ for any $i=1,...,n$, where $G_j$ is a
function of Riemann invariants satisfying
$$\partial_{r_i}G_j=\frac{\partial_{r_i}\phi_j}{\tan(\phi_i-\phi_j)}.$$
Here $v^{\bot}$ stands for the vector field rotated from $v$ by
$90^\circ$ counterclockwise.
\end{theorem}
We shall see in lemma below that the conditions (R) and ($\Phi$) are
almost equivalent. This lemma enables us to prove the following
theorem which is a generalization to our case of the result of
\cite{Sev}.
\begin{theorem}
Given any strictly Hyperbolic diagonal system
$$\cos \phi_i (r_i)_x+ \sin \phi_i (r_i)_y=0,\ i=1,...,n.$$
Then the condition ($\Phi$) is satisfied if and only if the system
can be written in the form of $n$ conservation laws
$$(g_i)_x+(h_i)_y=0, i=1,...,n.$$
%where the Jacobi matrices
%$\left(\partial_{r_j}{g_i}\right)$
%or $\left(\partial_{r_j}{h_i}\right)$ is non-degenerate.
\end{theorem}

\begin{remark}Saying that system (1) can be written in a
certain form, one means by definition that this form can be achieved
by a change of variables and multiplication by an invertible matrix:
Let $u=\Phi(w)$ be a regular change of variables and $C(w)$ is an
invertible matrix then one can write system (1) in the form
$$C(w)A(\Phi(w))D\Phi(w)w_x+C(w)B(\Phi(w))D\Phi(w)w_y=0,$$
where $D$ is the differential.
\end{remark}
%%%%%%%%%%%%%%%%%%%%%%%%%%%%%%%%%%%%%%%%%%%%%%%%%%%%%%%%%%%%%%%%%%%%%%%%%%
\section*{Acknowledgements}
It is a pleasure to thank Marshall Slemrod for his interest in the
results of this paper. I am also thankful to Andrey E. Mironov who
encouraged me not to be afraid of heavy computations.

%%%%%%%%%%%%%%%%
\section { Derivation along characteristics. Proof of Theorem 1.2. }
Differentiate the j-th equation of (2) with respect to the field
$v_j^\bot$. We have
$$0=L_{v_j^\bot}L_{v_j}r_j=L_{v_j}L_{v_j^\bot}r_j-L_{[v_j,v_j^\bot]}r_j.\eqno (3)$$
Compute now the derivative along the commutator:

$$
L_{[v_j,v_j^\bot]}r_j=L_{v_j}L_{v_j^\bot}r_j-L_{v_j^\bot}L_{v_j}r_j=$$
$$=L_{v_j}(-\sin\phi_j (r_j)_x+\cos\phi_j(r_j)_y)-L_{v_j^\bot}(\cos\phi_j
(r_j)_x+\sin\phi_j(r_j)_y)=$$
$$=(r_j)_x(-\cos^2
\phi_j(\phi_j)_x-\cos\phi_j\sin\phi_j(\phi_j)_y)+$$
$$+(r_j)_y(-\sin\phi_j\cos\phi_j(\phi_j)_x-\sin^2\phi_j(\phi_j)_y)+
$$
$$+(r_j)_x(-\sin^2
\phi_j(\phi_j)_x+\sin\phi_j\cos\phi_j(\phi_j)_y)+$$
$$+(r_j)_y(\sin\phi_j\cos\phi_j(\phi_j)_x-\cos^2\phi_j(\phi_j)_y)=
$$
$$=-(r_j)_x(\phi_j)_x-(r_j)_y(\phi_j)_y \eqno (4)$$

Notice that the derivatives $(r_j)_x, (r_j)_y$ can be expressed by
the following two identities:
$$\cos\phi_j
(r_j)_x+\sin\phi_j(r_j)_y=0,$$
$$-\sin\phi_j (r_j)_x+\cos\phi_j(r_j)_y=L_{v_j^\bot}r_j.$$
Therefore
$$(r_j)_x=-\sin\phi_jL_{v_j^\bot}r_j,$$
$$
(r_j)_x=\cos\phi_jL_{v_j^\bot}r_j. \eqno(5)$$ Substituting back to
(4) we get
$$L_{[v_j,v_j^\bot]}r_j=(L_{v_j^\bot}r_j)(\sin\phi_j
(\phi_j)_x-\cos\phi_j(\phi_j)_y)=-(L_{v_j^\bot}r_j)
(L_{v_j^\bot}\phi_j) \eqno(6)$$ By the chain rule for
$L_{v_j^\bot}\phi_j$ the last equation can be rewritten as follows
$$L_{[v_j,v_j^\bot]}r_j=-L_{v_j^\bot}r_j \sum_{i=1}^n
(\partial_{r_i}\phi_j)(L_{v_j^\bot}r_i)=$$
$$-(L_{v_j^\bot}r_j)^2(\partial_{r_j}
\phi_j)-(L_{v_j^\bot}r_j)\sum_{i\neq
j}(\partial_{r_i}\phi_j)(L_{v_j^\bot}r_i).\eqno(7)
$$
Let me express now the derivative $$L_{v_j^\bot}r_i=-\sin\phi_j
(r_i)_x+\cos\phi_j(r_i)_y, \eqno(8)$$ via $L_{v_j}r_i$ as follows.
Write
$$\cos\phi_i
(r_i)_x+\sin\phi_i(r_i)_y=0,$$
$$\cos\phi_j
(r_i)_x+\sin\phi_j(r_i)_y=L_{v_j}r_i. \eqno (9)$$
From these two
identities we have
$$(r_i)_x=\frac {\sin \phi_i}{\sin(\phi_i-\phi_j)}L_{v_j}r_i, \quad (r_j)_y=-\frac {\cos
\phi_i}{\sin(\phi_i-\phi_j)}L_{v_j}r_i. \eqno(10)
$$
Substitute expressions  (10) into (8) to get:
$$
L_{v_j^\bot}r_i=-\frac{L_{v_j}r_i}{tan(\phi_i-\phi_j)}. \eqno(11)
$$
Denote by
$$ w_i:=L_{v_i^\bot}r_i.
$$
Plug this together with (11) into equation (7) and then to (3):
$$
L_{v_j}(w_j)+(\partial_{r_j}\phi_j) (w_j)^2-w_j \sum_{i\neq
j}(\partial_{r_i}\phi_j)\frac{1}{tan(\phi_i-\phi_j)}L_{v_j}r_i=0.
\eqno(12)
$$
By Richness  ($\Phi$) we have that for all $j=1,...,n$ there exist
functions
$$G_j(r_1,...,r_n):\quad
\partial_{r_i}G_j=\frac{(\partial_{r_i}\phi_j)}{tan(\phi_i-\phi_j)},\
i\neq j.\eqno(13)
$$
By (13) we can rewrite (12) as
$$
L_{v_j}(w_j)+(\partial_{r_j}\phi_j) (w_j)^2-w_j \sum_{i\neq
j}(\partial_{r_i}G_j)L_{v_j}r_i=0,
$$
which is the same as
$$L_{v_j}(w_j)+(\partial_{r_j}\phi_j) (w_j)^2-w_j L_{v_j}G_j=0.
\eqno(14)
$$
Multiplying (14) by $exp(-G_j)$ we get the Riccati equation of the
first theorem for $j$ instead of $i$. This completes the proof of
theorem 1.2.
%%%%%%%%%%%%%%%%%%%%%%%%%%%%%%%%%%%%%%%%%%%%%%%%%%%%%%%%%%%%%%%%%%%%%%%%%%%%%%%%%%%%%%%%

\section{Conservation laws. Proof of theorem 1.3.}
We shall need the following key observation.
\begin{lemma}
Given two sets of functions $\lambda_i(r_1,...,r_n);\
\phi_i(r_1,...r_n),\ i=1,...,n)$ such that
$$\lambda_i\neq\lambda_j,\ \phi_i \neq\pi/2\ (mod \pi),\ \lambda_i=\tan\phi_i.$$
Then the conditions (R) and ($\Phi$) are equivalent.
\end{lemma}
The proof which I know is computational. It would be interesting to
find more conceptual proof.
\begin{proof}
Let us prove first that ($R$) implies ($\Phi$).

Denote by
$$a_{ij}:=
\frac{\partial_{r_i}\lambda_j}{\lambda_i-\lambda_j}.$$ Then $a_{ij}$
satisfy the following identities (\cite{Serre}):
$$
\partial_{r_i}a_{kj}=\partial_{r_k}a_{ij}=a_{ki}a_{ij}+a_{ik}a_{kj}-a_{kj}a_{ij}.
\eqno(15)
$$
In order to prove them differentiate with respect to $r_k$ the
identity

$\partial_{r_i}\lambda_j=a_{ij}(\lambda_i-\lambda_j)$ then
interchange the order of $i,k$, subtract one from the other and
divide by $\lambda_i-\lambda_k$.

Denote by
$$b_{ij}:=\frac{\partial_{r_i}\phi_j}{\tan(\phi_i-\phi_j)}=\frac{\partial_{r_i}\lambda_j}{\lambda_i-\lambda_j}
\frac{1+\lambda_i\lambda_j}{1+\lambda_j^2}=a_{ij}\frac{1+\lambda_i\lambda_j}{1+\lambda_j^2}.
\eqno (16)
$$
To prove ($\Phi$) we have to verify that the difference
$$d=\partial_{r_k}b_{ij}-\partial_{r_i}b_{kj}$$
vanishes. Let us compute $d$ explicitly:
$$d=(\partial_{r_k}a_{ij})\frac{1+\lambda_i\lambda_j}{1+\lambda_j^2}-(\partial_{r_i}a_{kj})\frac{1+\lambda_k\lambda_j}{1+\lambda_j^2}+
$$
$$+a_{ij}\partial_{r_k}\left(\frac{1+\lambda_i\lambda_j}{1+\lambda_j^2}\right)-a_{kj}\partial_{r_i}
\left(\frac{1+\lambda_k\lambda_j}{1+\lambda_j^2}\right).$$ By the
identities (15) and the condition (R) we have
$$d=(a_{ki}a_{ij}+a_{ik}a_{kj}-a_{kj}a_{ij})\frac{\lambda_j(\lambda_i-\lambda_k)}{1+\lambda_j^2}+$$

$$+a_{ij}\frac
{\partial_{r_k}(\lambda_i\lambda_j)(1+\lambda_j^2)-(1+\lambda_i\lambda_j)2\lambda_j\partial_{r_k}(\lambda_j)}{(1+\lambda_j^2)^2}-$$
$$
-a_{kj}\frac
{\partial_{r_i}(\lambda_k\lambda_j)(1+\lambda_j^2)-(1+\lambda_k\lambda_j)2\lambda_j\partial_{r_i}(\lambda_j)}{(1+\lambda_j^2)^2}.
\eqno(17)$$

Substitute now into the nominators of (17) the following expressions
for the derivatives of $\lambda_j$ from the definition of $a_{ij}$:
$$\partial_{r_i}\lambda_j=a_{ij}(\lambda_i-\lambda_j).$$
Then one has
$$d=(a_{ki}a_{ij}+a_{ik}a_{kj}-a_{kj}a_{ij})\frac{\lambda_j(\lambda_i-\lambda_k)}{1+\lambda_j^2}+$$
$$+a_{ij}\frac{(a_{kj}(\lambda_k-\lambda_j)\lambda_i+\lambda_j
a_{ki}(\lambda_k-\lambda_i)}{1+\lambda_j^2}-
a_{kj}\frac{(a_{ij}(\lambda_i-\lambda_j)\lambda_k+\lambda_j
a_{ik}(\lambda_i-\lambda_k)}{1+\lambda_j^2}-$$
$$-2a_{ij}\frac{(1+\lambda_i\lambda_j)\lambda_j
a_{kj}(\lambda_k-\lambda_j)}{(1+\lambda_j^2)^2}+ %%%%%%%%%%%%%%%
2a_{kj}\frac{(1+\lambda_k\lambda_j)\lambda_j
a_{ij}(\lambda_i-\lambda_j)}{(1+\lambda_j^2)^2}. \eqno(18)
$$
Notice that the identity (18) is a quadratic expression in
$a_{ij}s$. Collecting the coefficients of $a_{ij}a_{kj} ,\
a_{ik}a_{kj},\ a_{ki}a_{ij}$ one
comes to $d=0$. This proves lemma in one direction. %%%%%%%%%%%%%%%

Proof of the converse statement is very much analogous but with even
harder computations. I shall reproduce them sketchy. So we assume
the identities ($\Phi$) are satisfied. First one can obtain the
identity analogous to (15) for the derivatives
$\partial_{r_k}b_{ij}$ in the following way. Write
$$\partial_{r_i}\phi_j=b_{ij}\tan(\phi_i-\phi_j)=
b_{ij} \frac{\lambda_i-\lambda_j}{1+\lambda_i\lambda_j},\
\partial_{r_i}\lambda_j=b_{ij}\frac{(1+\lambda_j^2)(\lambda_i-\lambda_j)}{1+\lambda_i\lambda_j}.
\eqno(19).$$ Differentiating the first equality of (19) with respect
to $r_k$, using the identities (19) again and taking into account
(16) one has
$$\partial_{r_k}\partial_{r_i}\phi_j=\partial_{r_k}(b_{ij})\frac{(\lambda_i-\lambda_j)}{1+\lambda_i\lambda_j}+$$
$$+b_{ij}\left(1+\frac{(\lambda_i-\lambda_j)^2}{(1+\lambda_i\lambda_j)^2}\right)\left(b_{ki}
\frac{\lambda_k-\lambda_i}{1+\lambda_k\lambda_i}-b_{kj}
\frac{\lambda_k-\lambda_j}{1+\lambda_k\lambda_j}\right).
$$
Interchanging in this identity the order of indexes $i$ and $k$ and
using $\partial_{r_k}b_{ij}=\partial_{r_i}b_{kj}$ one has the
identity:
$$\partial_{r_k}(b_{ij})\frac{(\lambda_i-\lambda_k)(1+\lambda_j^2)}{(1+\lambda_i\lambda_j)(1+\lambda_k\lambda_j)}=$$
$$=b_{kj}\left(1+\frac{(\lambda_k-\lambda_j)^2}{(1+\lambda_k\lambda_j)^2}\right)\left(b_{ik}
\frac{\lambda_i-\lambda_k}{1+\lambda_i\lambda_k}-b_{ij}
\frac{\lambda_i-\lambda_j}{1+\lambda_i\lambda_j}\right)-$$
$$-b_{ij}\left(1+\frac{(\lambda_i-\lambda_j)^2}{(1+\lambda_i\lambda_j)^2}\right)\left(b_{ki}
\frac{\lambda_k-\lambda_i}{1+\lambda_i\lambda_k}-b_{kj}
\frac{\lambda_k-\lambda_j}{1+\lambda_k\lambda_j}\right).\eqno(20)$$
In order to verify (R) one computes
$$\partial_{r_k}a_{ij}-\partial_{r_i}a_{kj}=\partial_{r_k}b_{ij}
\frac{\lambda_j(\lambda_k-\lambda_i)(1+\lambda_j^2)}{(1+\lambda_i\lambda_j)(1+\lambda_k\lambda_j)}+$$
$$+b_{ij}\partial_{r_k}\left( \frac{1+\lambda_j^2}{1+\lambda_i\lambda_j}
\right)-b_{kj}\partial_{r_i}\left(
\frac{1+\lambda_j^2}{1+\lambda_k\lambda_j}\right).\eqno(21)
$$
The last step is to plug into (21) the expression (20) and also to
differentiate the last two brackets of (21) using the expression for
the derivatives $(19)$. Then one finally gets a quadratic expression
in $b_{ij}s$. Collecting similar terms one verifies that the right
hand side of (21) vanishes. Therefore (R) holds true. This proves
the lemma.
\end{proof}
%%%%%%%%%%%%%%%%%%%%%%%%%%%%%%%%%%%%%%%
It is easy now to prove Theorem 1.3.

\begin{proof}

Notice first of all that the statement of the second theorem is
local. Given a system which  is strictly Hyperbolic and is written
in the diagonal form
$$\cos \phi_i (r_i)_x+ \sin \phi_i (r_i)_y=0,\ i=1,...,n,$$
Let us give a proof first in one direction, namely assume that the
system can be written in the form of conservation laws
$$(g_i)_x+(h_i)_y=0, i=1,...,n.$$
Let me explain that then it must satisfy condition ($\Phi$). If
among $\phi_i$ there is one with $\cos\phi_i=0$ then one can apply a
small rotation of the plane $\mathbf{R}^2(x,y)$ and to get a new
system which has all angles different from $\pm\pi/2$. Notice that
the rotated system remains in the form of conservation laws and in
addition the differential $Dg$ becomes a non singular matrix, since
otherwise $\alpha=1,\beta=0$ would be the root of (P) but this is
impossible by $\phi_i\neq\pm\pi/2$. Denote
$$\lambda_i:=\tan\phi_i.$$ Use now Sevennec' theorem saying
that the diagonal system
$$(r_i)_x+\lambda_i (r_i)_y=0
$$ which can be written in the form of conservation laws
$$(g_i)_x+(h_i)_y=0, i=1,...,n$$
with the non-singular Jacobi matrix
$\left(\partial_{r_j}{g_i}\right)$ must satisfy (R). But by lemma in
this case (R) and ($\Phi$) are equivalent. So we get condition
($\Phi$) for rotated system. But this condition is obviously
rotationally invariant. Thus it holds also for the original system.

The proof in the opposite direction is analogous. First rotate the
plane exactly as above. Condition ($\Phi$)remains valid since it is
rotationally invariant. Then by the lemma (R) is valid as well and
then by Sevennec' theorem the rotated system can written in the form
of conservation laws. But then obviously the original one as well.
This completes the proof.

\end{proof}

%%%%%%%%%%%%%%%%%%%%%%%%%%%%%%%%%%%%%%%

\section{Geometric example.}
In this section we give a geometric example where the results of the
previous sections are important.

 Let $\rho$ be a Riemannian metric
on the 2-torus $\mathbb{T}^2=\mathbb{R}^2/\Gamma$, $\rho^t$ denotes
the geodesic flow. Assume that $\rho$ is written in conformal way:
$$ds^2=\Lambda(x,y)(dx^2+dy^2).$$
Let $F:T^* \mathbb{T}^2$ be a function on the cotangent bundle which
is homogeneous polynomial of degree $n$ with respect to the fibre:
$$ F=\sum_{k=0}^n
{a_{k}(x,y)}p_1^{n-k}p_2^{k}.
$$
 We are looking for such an $F$ which is an
integral of motion for the geodesic flow $\rho^t$, i.e.
$F\circ\rho^t=F$. We shall also assume that this $F$ is irreducible,
i.e of minimal possible degree. Let us mention that this problem is
classical there are very well studied examples of the geodesic flows
on the 2-torus which have integrals $F$ of degree one and two. We
refer to books \cite{BF} and \cite{Per} for the history and
discussion of this classical question with references therein. In
our recent papers with A.E.Mironov we used the so called
semi-geodesic coordinates. In these coordinates one arrives to a
remarkable Rich quasi-linear system of equations in evolution form
on the coefficients of the integral $F$ (\cite{BM1}, \cite{BM2}).
However it is very natural to be able to work in conformal
coordinates as well. In this case the quasi-linear system on the
coefficients has no evolution form any more but looks like:
$$A(U)U_x+B(U)U_y=0.$$
Let me write down explicitly the matrices for the case $n=3$ (this
case is already very interesting and not trivial see for example
\cite{DM}).

$$
 A(U)=  \left(
  \begin{array}{ccc}
   1 & 0 & 3a \\
   0 & 1 & 3b \\
  \Lambda & 0 & u\\
  \end{array}\right), \
  B(U)=\left(
  \begin{array}{ccc}
   0 & -1 & 3b \\
   1 & 0 & -3a \\
  0 & \Lambda & v\\
  \end{array}\right), \
  U=\left(
  \begin{array}{c}
   u\\
   v \\
  \Lambda \\
  \end{array}\right). \eqno (22)
$$
Here $a,b,u,v$ are related to the coefficients of the integral $a_i$
by the following
$$a_0=a+\frac{u}{\Lambda}, a_1=3b+\frac{v}{\Lambda}, a_2=-3a+\frac{u}{\Lambda}, a_3=-b+\frac{v}{\Lambda}.$$
It was noticed in \cite{Kol} that $a,b$ are in fact constants.
Computing polynomial $P=\det (\alpha B-\beta A)$ one has:
$$P=\alpha^3(v+3b\Lambda)+\alpha^2\beta(-u-9a\Lambda)+\alpha\beta^2(v-9b\Lambda)+\beta^3(-u+3a\Lambda).
$$
Let remark that it may happen at some points that both matrices
$A,B$ are degenerate, however polynomial $P$ for any point can not
vanish identically. This is because otherwise both constants $a,b$
vanish, but then one checks that in such a case the integral $F$ is
a product of the Hamiltonian with an integral of degree one in
momenta, therefore reducible.

Notice that quasi-linear system (22) is written in the form of
conservation laws $$(g_i)_x+(h_i)_y=0,$$
$$
g_1=u+3a\Lambda,g_2=v+3b\Lambda,g_3=u\Lambda,$$
$$h_1=-v+3b\Lambda,h_2=u-3a\lambda,
 h_3=v\Lambda.$$
 Moreover by a very general argument in
the Hyperbolic region this system can be written in the diagonal
form (2). Indeed introduce angular coordinate $\phi$ on the fibres
of the energy level $$\{\frac{1}{\Lambda}(p_1^2+p_2^2)=1\}:\
p_1=\sqrt{\Lambda}\cos\phi,p_2=\sqrt{\Lambda}\sin\phi,$$ then one
can verify that the condition on a function  $F$ to be an integral
of the flow reads
$$ F_x\cos\phi+F_y\sin\phi+F_\phi\left(\frac{\Lambda_y}{2\Lambda}\cos\phi-\frac{\Lambda_x}{2\Lambda}\sin\phi\right)=0$$
At the points where $F_\phi$ vanishes this equation takes particularly nice form:
$$ F_x\cos\phi+F_y\sin\phi=0.$$ Therefore critical values of $F$ on the fibre are Riemann invariants.
One can check also that the polynomial $P$ is proportional in fact
to the derivative of $F$ in the direction of the fibre. Moreover one
can check, as we did in \cite{BM1}, that in the Hyperbolic region
Riemann invariants form a regular change of variables. As a
consequence of Theorem 1.3 one concludes that in the Hyperbolic
region the system of this example is Rich in our generalized sense.
And therefore Theorem 1.2 tells us that the Riccati equation along
characteristics applies.

\section{Questions}
Several questions are very natural:

1. It would be interesting to find more conceptual proof of the
lemma in the framework differential-geometric approach by
Dubrovin-Novikov \cite{DN}.

2. How does generalized Hodograph method by Tsarev \cite{Tsarev}
work in our case?

3. How to analyze the behavior of the Riccati equation for the
example of previous section? It seems that genuine non-linearity
condition can not be expected for all eigenvalues.

\end{document}